\documentclass[fleqn,11pt]{article}
\usepackage{graphicx}
\usepackage{amscd,amssymb}
\usepackage{amsmath}
\title{The Mayer-Vietoris Property in Differential Cohomology}
\author{James Simons and Dennis Sullivan}
\date{}

\begin{document}

\maketitle

\begin{abstract}
In [1] it was shown that $K $ $\hat{}$, a certain differential cohomology
functor associated to complex K-theory, satisfies the Mayer-Vietoris
property when the underlying manifold is compact.  It turns out that
this result is quite general.  The work that follows shows the M-V
property to hold on compact manifolds for any differential cohomology
functor  $J$ $\hat{}$ associated to any  Z-graded cohomology functor  J( ,Z)
which, in each degree, assigns to a point a finitely generated group.
The approach is to show that the result follows from Diagram 1, the
commutative diagram we take as a definition of differential cohomology, and
Diagram 2, which combines the three Mayer-Vietoris
sequences  for $J^*( ,Z), J^*( ,R)$ and $J^*( ,R/Z)$.

\end{abstract}

Let $J = \sum \oplus J^{k}$ be a graded generalized cohomology functor.  We assume each $J^{k}(\textrm{point})$ is
finitely generated.  By a differential cohomology functor
associated to $J$ we mean a functor $\hat{J}$ on the category of
smooth manifolds with corners, together with four natural
transformations, $i_{1}, i_{2}, \delta_{1}, \delta_{2}$, which satisfies the following
commutative diagram of abelian group valued functors.

\underline{\textbf{Diagram 1}}
% \vspace{.5cm}
\begin{center}
\setlength{\unitlength}{0.5cm}
\begin{picture}(24,16)(2.5,0)
\thicklines
\put(5,1){$0$}

\put(20.5,1){$0$}

\put(6,2){\vector(1,1){1.5}}
\put(18,3.5){\vector(1,-1){1.5}}

\put(8,4.5){$\mathbf{\Lambda}^{k-1}/\mathbf{\Lambda}^{k-1}_{J}$}
\put(12,4.5){\vector(1,0){2.5}}
\put(16.5,4.5){$\mathbf{\Lambda}^{k}_{J}$}
\put(13,5){\small{$d$}}

\put(6.5,7.5){\vector(1,-1){1.5}}
\put(7.25,7){\small{$deRh$}}
\put(10.5,7){\small{$i_{2}$}}
\put(10.5,6){\vector(1,1){1.5}}
\put(14.5,7.5){\vector(1,-1){1.5}}
\put(15.5,7){\small{$\delta_{1}$}}
\put(17.75,7){\small{$deRh$}}
\put(18.5,6){\vector(1,1){1.5}}

\put(3,8){$\mathbf{H}^{k-1}(\cdot, \mathbb{R})$}
\put(12.75,8){$\hat{J}^{k}$}
\put(20,8){$\mathbf{H}^{k}(\cdot, \mathbb{R})$}

\put(6,10.5){\small{$p$}}
\put(6,9.5){\vector(1,1){1.5}}
\put(10.5,11){\vector(1,-1){1.5}}
\put(11.5,10.5){\small{$i_{1}$}}
\put(14,10.5){\small{$\delta_{2}$}}
\put(14.0,9.5){\vector(1,1){1.5}}
\put(18,11){\vector(1,-1){1.5}}
\put(19,10.5){\small{$i_{\mathbb{R}} \circ ch$}}

\put(7,12){$J^{k-1}(\cdot, \mathbb{R}/\mathbb{Z})$}
\put(12,12){\vector(1,0){2.5}}
\put(15.5,12){$J^{k}(\cdot, \mathbb{Z})$}
\put(13,12.5){\small{$b$}}

\put(5.5,14.5){\vector(1,-1){1.5}}
\put(18,13){\vector(1,1){1.5}}

\put(4.5,15){$0$}
\put(20,15){$0$}
\end{picture}
\end{center}

In the above the diagonals are short exact, and the upper and lower
four-term sequences are also exact, and
\begin{eqnarray*}
\mathbf{H}^{k}(\cdot \,, \mathbb{R}) & = & \sum_{j=0} \oplus H^{j}(\cdot \,,
J^{k-j}(\textrm{point}, \mathbb{R} ) \\
\\
\mathbf{\Lambda}^{k} & = & \sum_{j=0} \oplus \wedge^{j}(\cdot \,,
J^{k-j}(\textrm{point}, \mathbb{R} ) \\
\\
\mathbf{\Lambda}^{k}_{J} & = & (\textrm{de Rham})^{-1}(\textrm{Im}(ch \circ i_{\mathbb{R}}))
\end{eqnarray*}
$ch : J^{k}(\mathbb{Z}) \rightarrow \mathbf{H}^{k}(\mathbb{Q})$ is the
canonical map, $i_{\mathbb{R}}$ is induced by $\mathbb{Q} \rightarrow \mathbb{R}$, $p$ is induced by the coefficient
sequence $\mathbb{Z} \rightarrow \mathbb{R} \rightarrow
\mathbb{R}/\mathbb{Z}$, and $b$ denotes the associated Bockstein map.
The maps $deRh$ and $d$ are respectively the de Rham map and the
exterior differential.

\underline{\textbf{Theorem (Mayer-Vietoris Property)}}:  Let $X$ be a
compact smooth manifold.
Assume $X = A \cup B$, and $A
\cap B = D$, a co-dim $0$ submanifold with collar neighborhoods in
both $A$ and $B$.  Then, if $f_{A} \in \hat{J}^{k}(A)$ and $f_{B} \in \hat{J}^{k}(B)$
with $f_{A} | D = f_{B} | D$, then there exists $f \in
\hat{J}^{k}(X)$ with $f | A = f_{A}$ and $f | B = f_{B}$.

\underline{\textbf{Proof}}:  Since $\delta_{2}(f_{A}) | D =
\delta_{2}(f_{B}) | D$, the Mayer-Vietoris property for $J$ shows there
exists $v \in J^{k}(X)$ with $v | A = \delta_{2}(f_{A})$ and $v | B =
\delta_{2}(f_{B})$.  Choose $h \in \hat{J}^{k}(X)$ with $\delta_{2}(h) =
v$.  By naturality
\begin{eqnarray*}
\delta_{2}(h | A) = \delta_{2}(f_{A}) & \textrm{and} & \delta_{2}(h |
B) = \delta_{2}(f_{B})
\end{eqnarray*}
thus by Diagram 1

1) \hspace*{.3cm} $h | A - f_{A} = i_{2} (\{\alpha_{A}\})$
\hspace{1cm} $\{\alpha_{A}\} \in \mathbf{\Lambda}^{k-1}(A)/\mathbf{\Lambda}^{k-1}_{J}(A)$ \\

\vspace*{-.5cm}
\hspace*{.7cm} $ h | B - f_{B} = i_{2} (\{\alpha_{B}\})$ \hspace{1cm} $\{\alpha_{B}\} \in \mathbf{\Lambda}^{k-1}(B)/\mathbf{\Lambda}^{k-1}_{J}(B)$.

Under restriction to $D$ the left hand sides are equal by hypothesis.
Since $i_{2}$ is an injection, by naturality we see that
$\{\alpha_{A}\} | D = \{\alpha_{B}\} | D$.

Suppose one can find $\bar{h} \in \hat{J}^{k}(X)$ with $\bar{h} | A =
i_{2}(\{\alpha_{A}\})$ and $\bar{h} | B = i_{2}(\{\alpha_{B}\})$.
Then by 1), \, ($h - \bar{h}) | A = f_{A}$ and $(h - \bar{h}) | B = f_{B}$,
and the problem is solved.

The problem thus reduces to the case that $f_{A} =
i_{2}(\{\alpha_{A}\})$ and $f_{B} = i_{2}(\{\alpha_{B}\})$.  The
remainder of the proof will be restricted to this case.

Since $\{\alpha_{A}\} | D - \{\alpha_{B}\} | D = 0$ we must have
$\alpha_{A} | D - \alpha_{B} | D \in \mathbf{\Lambda}^{k-1}_{J}(D)$.  If instead
we had chosen $\alpha_{A}^{\prime}$ and $\alpha_{B}^{\prime}$
representing $\{\alpha_{A}\}$ and $\{\alpha_{B}\}$ then $\alpha_{A} -
\alpha_{A}^{\prime} \in \mathbf{\Lambda}^{k-1}_{J}(A)$ and $\alpha_{B} -
\alpha_{B}^{\prime} \in \mathbf{\Lambda}^{k-1}_{J}(B)$.  Thus
\begin{eqnarray*}
(\alpha_{A} | D - \alpha_{B} | D) - (\alpha_{A}^{\prime} | D -
\alpha_{B}^{\prime} | D)  \in \mathbf{\Lambda}^{k-1}_{J}(A) | D + \mathbf{\Lambda}^{k-1}_{J}(B) | D
\end{eqnarray*}
and therefore
\begin{eqnarray*}
w(\{\alpha_{A}\},\{\alpha_{B}\}) = (\alpha_{A} | D - \alpha_{B} | D)
\in \frac{\mathbf{\Lambda}^{k-1}_{J}(D)}{\mathbf{\Lambda}^{k-1}_{J}(A) | D + \mathbf{\Lambda}^{k-1}_{J}(B) | D }
\end{eqnarray*}
is well defined.  Suppose $w(\{\alpha_{A}\},\{\alpha_{B}\}) = 0$.
Then $\alpha_{A} | D - \alpha_{B} | D = \beta_{A} | D - \beta_{B} |
D$, where $\beta_{A} \in \mathbf{\Lambda}^{k-1}_{J}(A)$ and $\beta_{B} \in
\mathbf{\Lambda}^{k-1}_{J}(B)$.  Thus $\{\alpha_{A}\} = \{\alpha_{A} -
\beta_{A}\}$, $\{\alpha_{B}\} = \{\alpha_{B} - \beta_{B}\}$, and
\begin{eqnarray*}
(\alpha_{A} - \beta_{A}) | D = (\alpha_{B} - \beta_{B}) | D.
\end{eqnarray*}
Since $D$ has co-dim $0$ and collar neighborhoods in both $A$ and $B$,
there exists a unique $\theta \in \mathbf{\Lambda}^{k-1}(X)$ with $\theta | A =
\alpha_{A} - \beta_{A}$ and $\theta | B = \alpha_{B} - \beta_{B}$.
Thus $\{\theta\} | A = \{\alpha_{A}\}$ and $\{\theta\} | B =
\{\alpha_{B}\}$, which implies that $i_{2}(\{\theta\}) | A =
i_{2}(\{\alpha_A\})$ and $i_{2}(\{\theta\}) | B = i_{2}(\{\alpha_B\})$.

We have therefore shown

2) \hspace*{.3cm} $w(\{\alpha_A\},\{\alpha_B\}) = 0  \, \implies \, \textrm{problem is solved.}$

Set $J^{k}_{o}(X) = \{ v \in J^{k}(X) \hspace{.25cm} | \hspace{.25cm} v|A = 0 = v|B\}$.  Let $v \in
J^{k}_{o}(X)$ and choose $h \in \hat{J}^{k}(X)$ with $\delta_{2}(h) =
v$.  By naturality, $\delta_{2}(h|A) = 0 = \delta_{2}(h|B)$.  Thus
\begin{eqnarray*}
h|A & = & i_{2}(\{\gamma_A\}), \\
h|B & = & i_{2}(\{\gamma_B\}), \, \textrm{and} \\
\{\gamma_A\} | D & = & \{\gamma_B\} | D.
\end{eqnarray*}
Set
\[
\Omega(v) = w(\{\gamma_A\}, \{\gamma_B\}).
\]
To see that $\Omega$ is well defined, let $\bar{h} \in \hat{J}^{k}(X)$
with $\delta_{2}(\bar{h}) = v$.  Then $\bar{h} = h + i_{2}(\{\rho\})$
for some $\rho \in \mathbf{\Lambda}^{k-1}(X)$.  So $\bar{h} | A = i_{2} \{
\gamma_{A} + \rho | A\}$ and $\bar{h} | B = i_{2} \{\gamma_{B} + \rho | B\}$.
Since $(\rho | A) | D = \rho |D = (\rho | B) | D$, the definition of $w$
shows $w(\{\gamma_{A} + \rho | A \}, \{\gamma_{B} + \rho | B \} =
w(\{\gamma_{A}\},\{\gamma_{B}\})$.  Thus,
\[
\Omega: \, {J}^{k}_{o}(X) \rightarrow \frac{\mathbf{\Lambda}^{k-1}_{J}(D)}{\mathbf{\Lambda}^{k-1}_{J}(A) | D + \mathbf{\Lambda}^{k-1}_{J}(B) | D }
\]
is well defined, and is clearly a homomorphism.

Now, given $\{\alpha_{A}\},\{\alpha_{B}\}$ with $\{\alpha_{A}\} | D =
\{\alpha_{B}\} | D$, suppose we can find $v \in J^{k}_{o}(X)$ with
$\Omega(v) = w(\{\alpha_{A}\},\{\alpha_{B}\})$.  Then, choosing $h$
with $\delta_{2}(h) = B$ and letting $h | A = \{\gamma_{A}\}$ and $h | B
= \{\gamma_{B}\}$, we see
\[
w(\{\alpha_{A} - \gamma_{A}\},\{\alpha_{B} - \gamma_{B}\}) = 0.
\]
By 2), this implies there exists $\theta \in \mathbf{\Lambda}^{k-1}(X)$ with
\begin{eqnarray*}
\{\theta\} | A & = & \{\alpha_{A} - \gamma_{A}\} = \{\alpha_{A}\} -
\{\gamma_{A}\} \\
\{\theta\} | B & = & \{\alpha_{B} - \gamma_{B}\} = \{\alpha_{B}\} -
\{\gamma_{B}\}
\end{eqnarray*}
Thus
\begin{eqnarray*}
(i_{2}(\{\theta\}) + h) | A = i_{2}(\{\alpha_{A}\}) \\
(i_{2}(\{\theta\}) + h) | B = i_{2}(\{\alpha_{B}\})
\end{eqnarray*}
and so $i_{2}(\{\theta\}) + h$ solves the problem for the coherent
pair $i_{2}(\{\alpha_{A}\}), i_{2}(\{\alpha_{B}\})$.

The proof of the theorem will clearly be completed if we can show

$\ast)$ \hspace{0.3cm} $\Omega$ is surjective.

The remainder of the work will be devoted to proving $\ast$).

We consider the following diagram in which the rows are Mayer-Vietoris
exact sequences of the various cohomology functors.

\underline{\textbf{Diagram 2}}
\[ \begin{array}{ccccccc}
\scriptstyle J^{k-2}(A,\mathbb{R}/\mathbb{Z}) \oplus
 J^{k-2}(B,\mathbb{R}/\mathbb{Z}) & \xrightarrow{\Delta_{1}} &
\scriptstyle J^{k-2}(D,\mathbb{R}/\mathbb{Z}) & \xrightarrow{d_{1}^{\ast}} &
\scriptstyle J^{k-1}(X,\mathbb{R}/\mathbb{Z}) & \xrightarrow{\sum_{1}} &
\scriptstyle J^{k-1}(A,\mathbb{R}/\mathbb{Z}) \oplus
 J^{k-1}(B,\mathbb{R}/\mathbb{Z}) \\

\downarrow \scriptstyle b & & \downarrow \scriptstyle b & & \downarrow
\scriptstyle b & & \downarrow \scriptstyle b \\

\scriptstyle J^{k-1}(A,\mathbb{Z}) \oplus J^{k-1}(B,\mathbb{Z}) & \xrightarrow{\Delta_{2}} &
\scriptstyle J^{k-1}(D,\mathbb{Z}) & \xrightarrow{d_{2}^{\ast}} &
\scriptstyle J^{k}(X,\mathbb{Z}) & \xrightarrow{\sum_{2}} &
\scriptstyle J^{k}(A,\mathbb{Z}) \oplus J^{k}(B,\mathbb{Z}) \\

\downarrow \scriptstyle ch & & \downarrow \scriptstyle ch & &
\downarrow \scriptstyle ch & & \downarrow \scriptstyle ch \\

\scriptstyle \mathbf{H}^{k-1}(A,\mathbb{Q}) \oplus \mathbf{H}^{k-1}(B,\mathbb{Q}) & \xrightarrow{\Delta_{3}} &
\scriptstyle \mathbf{H}^{k-1}(D,\mathbb{Q}) & \xrightarrow{d_{3}^{\ast}} &
\scriptstyle \mathbf{H}^{k}(X,\mathbb{Q}) & \xrightarrow{\sum_{3}} &
\scriptstyle \mathbf{H}^{k}(A,\mathbb{Q}) \oplus \mathbf{H}^{k}(B,\mathbb{Q}) \\
\end{array} \]

The $\Delta$'s are the differences of the restrictions to $D$ of the
individual components.  $d^{\ast}$ is the Mayer-Vietoris promotion map.
$\sum$ restricts an element to each of $A$ and $B$ and takes their
direct sum.  $b$ is the Bockstein map, and $ch$ is defined in Diagram
1.  It is well known that all $2 \times 2$ boxes commute up to appropriate
sign in the graded sense.  Note that $\textrm{Im}(ch)$ is a spanning lattice in $\mathbf{H}^{\ast}(\cdot,
\mathbb{Q})$.

The proof of $\ast$) will now follow from a series of lemmas.

\vspace{0.3cm}

\underline{\textbf{Lemma 1}}:
\[
\frac{\mathbf{\Lambda}^{k-1}_{J}(D)}{\mathbf{\Lambda}^{k-1}_{J}(A) | D +
\mathbf{\Lambda}^{k-1}_{J}(B) | D }
\xrightarrow[\cong]{\textrm{\tiny{de Rham}}}
\frac{i_{\mathbb{R}} \circ
ch(J^{k-1}(D,\mathbb{Z}))}{i_{\mathbb{R}} \circ
ch(\textrm{Im}(\Delta_{2}))} \xleftarrow[\cong]{i_{\mathbb{R}}}
\frac{ch(J^{k-1}(D,\mathbb{Z}))}{ch(\textrm{Im}(\Delta_{2}))}
\]

\underline{\textbf{Proof}}:  Since de Rham maps the denominator of the
first expression into that of the second, the first map is well
defined and is onto since the map of the numerator is onto.  If
$\theta \in \mathbf{\Lambda}^{k-1}(D)$ maps to an element of $i_{\mathbb{R}}
\circ ch(\textrm{Im}(\Delta_{2}))$ there must be an $\eta \in
\mathbf{\Lambda}^{k-1}_{J}(A) | D + \mathbf{\Lambda}^{k-1}_{J}(B) | D$ and $\mu \in
\mathbf{\Lambda}^{k-1}_{exact}$ with $\theta = \eta + \mu$.  But any exact form
on $D$ is the restriction of an exact form on $A$.  Moreover
$\mathbf{\Lambda}^{\ast}_{exact} \subseteq \mathbf{\Lambda}^{\ast}_{J}$, and thus $\mu \in
\mathbf{\Lambda}^{k-1}_{J}(A) | D$, and so $\theta \in \mathbf{\Lambda}^{k-1}_{J}(A) | D +
\mathbf{\Lambda}^{k-1}_{J}(B) | D.$  Therefore de Rham is $1:1$, and so an
isomorphism.  That $i_{\mathbb{R}}$ induces an isomorphism is straightforward.
 \hspace{1 cm}   $\blacksquare$

By the above, we may consider

3) \hspace*{.3cm} $\displaystyle \Omega: \, J^{k}_{o}(X) \rightarrow \frac{ch(J^{k-1}(D,\mathbb{Z}))}{ch(\textrm{Im}(\Delta_{2}))}$.

Note that by Diagram 2,
\[
ch(\textrm{Im}(\Delta_{2})) = \Delta_{3}(\textrm{Im}(ch)) \subseteq \textrm{Im}(\Delta_{3}).
\]
\vspace{0.3cm}

\underline{\textbf{Lemma 2}}:  Let
\[
\phi: \,
\frac{ch(J^{k-1}(D,\mathbb{Z}))}{ch(\textrm{Im}(\Delta_{2}))}
\rightarrow \frac{\mathbf{H}^{k-1}(D,\mathbb{Q})}{\textrm{Im}(\Delta_{3})}
\]
be the map induced by inclusion.  Then
\[
\textrm{ker}(\phi) = \textrm{tor}\left(\frac{ch(J^{k-1}(D,\mathbb{Z}))}{ch(\textrm{Im}(\Delta_{2}))}\right).
\]

\underline{\textbf{Proof}}:  Clearly $\textrm{torsion} \subseteq
\textrm{ker}(\phi)$ since the image of $\phi$ lies in a rational
vector space.  Let $x_{1}, \cdots, x_{n}$ be a set of generators of
$J^{k-1}(A,\mathbb{Z}) \oplus J^{k-1}(B,\mathbb{Z})$.  Then $\{
ch(x_{i})\}$ span $\mathbf{H}^{k-1}(A,\mathbb{Q}) \oplus
\mathbf{H}^{k-1}(B,\mathbb{Q})$, and thus $\{\Delta_{3}(ch(x_{i}))\}$ =
$\{ch(\Delta_{2}(x_{i}))\}$ generate $\textrm{Im}(\Delta_{3})$.
Therefore if $y \in ch(J^{k-1}(D,\mathbb{Z}))$, and $y \in
\textrm{Im}(\Delta_{3})$, $y = \sum q_{i} ch(\Delta_{2}(x_{i}))$ for
some choice of rational $\{q_{i}\}$.  Clearing denominators leads to
integers $m, m_{1}, \cdots, m_{n}$ with $my = \sum m_{i}
ch(\Delta_{2}(x_{i}))$.  Thus $y$ represents a torsion element in
$ch(J^{k-1}(D,\mathbb{Z})) / ch(\textrm{Im}(\Delta_{2}))$.
 \hspace{1 cm} $\blacksquare$

Let
\[
\Pi: \, \frac{ch(J^{k-1}(D,\mathbb{Z}))}{ch(\textrm{Im}(\Delta_{2}))}
\rightarrow
\frac{ch(J^{k-1}(D,\mathbb{Z}))}{ch(\textrm{Im}(\Delta_{2}))} \, / \, \textrm{Torsion}
\]

\vspace{0.3cm}
\underline{\textbf{Lemma 3}}:  $\Pi \circ \Omega$ is surjective.

\underline{\textbf{Proof}}:   From Diagram 2 we derive
\begin{center}
\setlength{\unitlength}{0.5cm}
\begin{picture}(24,10)
\thicklines
\put(0,0.25){$\displaystyle \frac{\mathbf{H}^{k-1}(D,\mathbb{Q})}{\textrm{Im}(\Delta_{3})}$}
\put(8,0.75){$\scriptstyle d^{\ast}_{3}$}
\put(7,0.5){\vector(1,0){2.5}}
\put(8,0){$\scriptstyle 1:1$}
\put(11,0.5){$\mathbf{H}^{k}(X,\mathbb{Q})$}

\put(1.5,2.5){$\scriptstyle \phi$}
\put(2,3.25){\vector(0,-1){1.5}}
\put(2.25,2.5){$\scriptstyle 1:1$}

\put(0,4.75){$\displaystyle
\frac{ch(J^{k-1}(D,\mathbb{Z}))}{ch(\textrm{Im}(\Delta_{2}))} / \textrm{Torsion}$}
\put(12.25,4.75){$\scriptstyle ch$}
\put(13,6){\vector(0,-1){2.5}}

\put(0.5,7.25){$\scriptstyle \Pi \circ \frac{ch}{ch}$}
\put(2,8){\vector(0,-1){1.5}}
\put(2.25,7.25){\tiny{onto}}
\put(8.5,7.5){$\scriptstyle \Pi \circ \Omega$}
\put(11,8.5){\vector(-1,-1){2.5}}

\put(0,9){$\displaystyle \frac{J^{k-1}(D,\mathbb{Q})}{\textrm{Im}(\Delta_{2})}$}
\put(8,9.5){$\scriptstyle d^{\ast}_{2}$}
\put(7,9.25){\vector(1,0){2.5}}
\put(8,8.75){$\scriptstyle \cong$}
\put(11,9.25){$J^{k}_{o}(X,\mathbb{Z})$}

\end{picture}
\end{center}
By $\frac{ch}{ch}$ we mean the application of $ch$ to both numerator
and denominator.  Clearly $\Pi \circ \frac{ch}{ch}$ is onto.  Since
 $J^{k}_{o}(X,\mathbb{Z}) = \textrm{ker}(\sum_{2})$, Diagram 2 shows
that, as used above, $d^{\ast}_{2}$ is an isomorphism and
$d^{\ast}_{3}$ is $1:1$.  By Lemma 2, $\phi$ is $1:1$.

Recalling the definition of $w$ and $\Omega$, an element $v \in J^{k}_{o}(X,\mathbb{Z})$ and a choice of $h \in \hat{J}^{k}(X)$ with
$\delta_{1}(h) = v$ gives rise to elements $\gamma_{A}, \gamma_{B} \in
\mathbf{\Lambda}^{k-1}(A), \mathbf{\Lambda}^{k-1}(B)$, with $ \gamma_{A} | D - \gamma_{B} | D \in \mathbf{\Lambda}_{J}^{k-1}(D)$, the \mbox{de Rham} image of which lies
in $H^{k-1}(D,\mathbb{Q})$.  Letting $ [\gamma_{A} | D - \gamma_{B} |
D ]$ represent its rational cohomology class, and using 3), we may write
\[
\Omega(v) = [\gamma_{A} | D - \gamma_{B} | D] \, \bmod ch(\textrm{Im}(\Delta_{2})).
\]
From Diagram 1, we see $\delta_{1}(h | A) = d \gamma_{A}$, and
$\delta_{1}(h | B) = d \gamma_{B}$.  Since
\[
d \gamma_{A} | D - d \gamma_{B} | D = d(\gamma_{A} | D - \gamma_{B} |
D) = 0
\]
we may define the closed form $\eta$ on $X$ by $\eta | A =
d\gamma_{A}$ and $\eta | B = d\gamma_{B}$.  Clearly $\eta =
\delta_{1}(h)$ and thus $[\eta] = i_{\mathbb{R}} (ch(v))$.  Let $d^{\ast}$
denote the Mayer-Vietoris promotion map in $\mathbf{H}^{\ast}(\cdot,
\mathbb{R})$.  From the definition of $[\eta]$, we see that
\[
d^{\ast}(i_{\mathbb{R}}([\gamma_{A} | D - \gamma_{B} | D])) = i_{\mathbb{R}} (ch(v)).
\]
Since $d^{\ast} \circ i_{\mathbb{R}} = i_{\mathbb{R}} \circ
d_{3}^{\ast}$ we see that $d_{3}^{\ast}([\gamma_{A} | D - \gamma_{B} |
D]) = ch(v)$, and this implies that in the above diagram

4) \hspace*{.3cm} $ d_{3}^{\ast} \circ \phi \circ \Pi \circ \Omega(v) = ch(v)$.

To show that $\Pi \circ \Omega$ is surjective, let $x \in
ch(J^{k-1}(D,\mathbb{Z}))/ch(\textrm{Im}(\Delta_{2}))  \, \bmod \,
\textrm{torsion}$, and choose $y$ with $\Pi \circ \frac{ch}{ch}(y) =
x$.  Then by 4)
\[
d_{3}^{\ast} \circ \phi \circ \Pi \circ \Omega(d_{2}^{\ast}(y)) =
ch(d_{2}^{\ast}(y)) = d_{3}^{\ast} \circ ch(y) = d_{3}^{\ast} \circ
\phi \circ \Pi \circ \frac{ch}{ch}(y).
\]
Since $d_{3}^{\ast} \circ \phi$ is $1:1$ we must have $\Pi \circ
\Omega(d_{2}^{\ast}(y)) = \Pi \circ \frac{ch}{ch}(y) = x$.   \hspace{1 cm} $\blacksquare$

By the commutativity of Diagram 2 we note that
$b(\textrm{Im}(d^{\ast}_{1})) = d^{\ast}_{2}(\textrm{Im}(b)) \subseteq
\textrm{ker}(\sum_{2}) = J^{k}_{o}(X,\mathbb{Z})$.

\vspace{0.3cm}

\underline{\textbf{Lemma 4}}:  $b(\textrm{Im}(d^{\ast}_{1}))
\subseteq \textrm{ker}(\Omega)$.

\underline{\textbf{Proof}}:  Let $x \in
J^{k-2}(D,\mathbb{R}/\mathbb{Z})$, and $b(d^{\ast}_{1}(x)) =
v \in J^{k}_{o}(X,\mathbb{Z})$.  To compute $\Omega(v)$ we need $h \in
\hat{J}^{k}(X)$ with $\delta_{2}(h) = v$ and consider $h | A$ and $h |
B$.  By Diagram 1 we may take $h = i_{1}(d^{\ast}_{1}(x))$.  But
$i_{1}(d^{\ast}_{1}(x)) | A = i_{1}((d^{\ast}_{1}(x) | A) = 0$ since
$d^{\ast}_{1}(x) \in \textrm{ker}(\sum_{1})$.  Similarly for $B$.
Thus $\Omega(v) = 0$.   \hspace{1 cm} $\blacksquare$

Thus, we may regard

5) \hspace*{.3cm} $\displaystyle \Omega: \,
\frac{J^{k}_{o}(X,\mathbb{Z})}{b(\textrm{Im}(d^{\ast}_{1}))}
\rightarrow \frac{ch(J^{k-1}(D,\mathbb{Z}))}{ch(\textrm{Im}(\Delta_{2}))}$.

% AKH
%\newpage
\underline{\textbf{Lemma 5}}:
%\vspace{.5cm}
\begin{center}
\setlength{\unitlength}{0.5cm}
\begin{picture}(24,8)
%\begin{picture}(16,8)(2.5,0)
\thicklines
\put(12.5,7.25){$\displaystyle\frac{J^{k}(X,\mathbb{Z})}{b(\textrm{Im}(d^{\ast}_{1}))}$}

\put(14,5.75){\vector(0,-1){1.5}}
\put(14,4.75){\vector(0,1){1.5}}
\put(14.25,5){$\scriptstyle \cong$}
\put(12,2){$\displaystyle \frac{ch(J^{k-1}(D,\mathbb{Z}))}{ch(\textrm{Im}(\Delta_{2}))}$}

\put(10.75,4.0){$\scriptstyle ch$}
\put(10,4.5){\vector(1,-1){1.5}}
\put(10,3.25){$\scriptstyle \cong$}

\put(0,5){$\displaystyle \frac{J^{k-1}(D,\mathbb{Z})}{\textrm{Im}(\Delta_{2})
+ \textrm{Tor}(J^{k-1}(D,\mathbb{Z}))}$}

\put(10,6.5){$\scriptstyle d^{\ast}_{2}$}
\put(10,5.5){\vector(1,1){1.5}}
\put(10.75,5.75){$\scriptstyle \cong$}

\end{picture}
\end{center}
\underline{\textbf{Proof}}:  In the upper case we note that
$d^{\ast}_{2} : J^{k-1}(D,\mathbb{Z})/\textrm{Im}(\Delta_{2})
\xrightarrow{\cong} J^{k}_{o}(X,\mathbb{Z})$, and
$d^{\ast}_{2}(\textrm{Tor}(J^{k-1}(D,\mathbb{Z}))) =
d^{\ast}_{2}(\textrm{Im}(b)) = b(\textrm{Im}(d^{\ast}_{1}))$.  In the
lower case we note that $ch: \, J^{k-1}(D,\mathbb{Z}) /
\textrm{Tor}(J^{k-1}(D,\mathbb{Z})) \xrightarrow{\cong}
ch(J^{k-1}(D,\mathbb{Z}))$.  The vertical isomorphism then
follows.  \hspace{1 cm}  $\blacksquare$

\vspace{0.3cm}

\underline{\textbf{Lemma 6}}: \hspace*{.3cm}  $ \displaystyle \Omega: \,
\textrm{Tor}\left(\frac{J^{k-1}(D,\mathbb{Z})}{b(\textrm{Im}(d^{\ast}_{1}))}\right)
\xrightarrow{\cong}
\textrm{Tor}\left(\frac{ch(J^{k-1}(D,\mathbb{Z}))}{ch(\textrm{Im}(d^{\ast}_{1}))}\right)
$.

\underline{\textbf{Proof}}:  Since $b(\textrm{Im}(d^{\ast}_{1}))
\subseteq \textrm{Tor}(J^{k}_{o}(X,\mathbb{Z}))$,
\[
\textrm{Tor}\left(\frac{J^{k}_{o}(X,\mathbb{Z})}{b(\textrm{Im}(d^{\ast}_{1}))}\right)
=
\frac{\textrm{Tor}(J^{k}_{o}(X,\mathbb{Z}))}{b(\textrm{Im}(d^{\ast}_{1}))}
= \frac{\textrm{Im}(b) \cap \textrm{ker}(\sum_{2})}{b(\textrm{Im}(d^{\ast}_{1}))}
\]
Let $x \in \textrm{Im}(b) \cap \textrm{ker}(\sum_{2})$, i.e. $x =
b(u)$, where $u \in J^{k-1}(X,\mathbb{R}/\mathbb{Z})$ and $b(u) | A =
0 = b(u) | B$.  From Diagram 1 and naturality we see
$\delta_{2}(i_{1}(u) | A) = 0 = \delta_{2}(i_{1}(u) | B)$, and so
$i_{1}(u) | A = i_{2}(\{\theta_{A}\})$ and $i_{1}(u) | B =
i_{2}(\{\theta_{B}\})$.  Since $\delta_{1} \circ i_{1} = 0$ and
$\delta_{1} \circ i_{2} = d$, we see
\[
\theta_{A} \in \mathbf{\Lambda}^{k-1}_{\textrm{closed}}(A), \, \theta_{B} \in
\mathbf{\Lambda}^{k-1}_{\textrm{closed}}(B) \hspace{1cm} \textrm{and} \hspace{1cm}  \theta_{A} | D -
\theta_{B} | D \in \mathbf{\Lambda}^{k-1}_{J}(D).
\]
Using the original formulation of $\Omega$
\[
\Omega(x) = \theta_{A} | D - \theta_{B} | D \, \bmod
\mathbf{\Lambda}^{k-1}_{J}(A) | D + \mathbf{\Lambda}^{k-1}_{J}(B) | D.
\]
Now suppose $\Omega(x) = 0$.  This implies one can find $\gamma_{A}
\in \mathbf{\Lambda}^{k-1}_{J}(A)$ and $\gamma_{B} \in \mathbf{\Lambda}^{k-1}_{J}(B)$ with
\[
\theta_{A} | D - \theta_{B} | D = \gamma_{A} | D - \gamma_{B} | D.
\]
Since $\{\theta_{A}\} = \{\theta_{A} - \gamma_{A}\}$ and $\{\theta_{B}\} = \{\theta_{B} - \gamma_{B}\}$
we see
\begin{eqnarray*}
i_{1}(u) | A = i_{2}(\{\theta_{A} - \gamma_{A}\}) \\
i_{1}(u) | B = i_{2}(\{\theta_{B} - \gamma_{B}\})
\end{eqnarray*}
where $(\theta_{A} - \gamma_{A}) | D = (\theta_{B} - \gamma_{B}) |
D)$.  Thus we may define $\sigma \in
\mathbf{\Lambda}^{k-1}_{\textrm{closed}}(X)$ by $\sigma | A = \theta_{A} -
\gamma_{A}$ and $\sigma | B = \theta_{B} - \gamma_{B}$.

Let $[\sigma] \in \mathbf{H}^{k-1}(X,\mathbb{R})$ be the de Rham class
represented by $\sigma$.  Referring to Diagram 1 we have $p([\sigma])
\in J^{k-1}(X,\mathbb{R}/\mathbb{Z})$ and
\begin{eqnarray*}
i_{1}(u|A) & = & i_{2}(\textrm{deRh}([\sigma]|A))  = i_{1}(p([\sigma]|A)) \\
i_{1}(u|B) & = & i_{2}(\textrm{deRh}([\sigma]|B))  = i_{1}(p([\sigma]|B)).
\end{eqnarray*}
By injectivity of $i_{1}$ we see
\[
(u - p([\sigma])) | A = 0 = (u - p([\sigma])) | B.
\]
Thus $u - p([\sigma]) \in \textrm{Im}(d^{\ast}_{1})$.  Since
$\textrm{Im}(p) = \textrm{ker}(b)$ we see
\[
x = b(u) = b(u - p([\sigma])) \in b(\textrm{Im}(d^{\ast}_{1})).
\]
Thus $\Omega |
\textrm{Tor}(\frac{J^{k}_{o}(X,\mathbb{Z})}{b(\textrm{Im}(d^{\ast}_{1}))})$
is $1:1$.  By the assumption that $J^{k}(\textrm{point})$ is finitely
generated and $X$ is compact it follows that
$\textrm{Tor}\left(\frac{J^{k}_{o}(X,\mathbb{Z})}{b(\textrm{Im}(d^{\ast}_{1}))}\right)$
is finite.

By Lemma 5
\[
\textrm{card}\left(\textrm{Tor}\left(\frac{J^{k}_{o}(X,\mathbb{Z})}{b(\textrm{Im}(d^{\ast}_{1}))}\right)\right)
= \textrm{card}\left(\textrm{Tor}\left(\frac{ch(J^{k-1}(D,\mathbb{Z}))}{ch(\textrm{Im}(\Delta_{2}))}\right)\right).
\]
Surjectivity thus follows from injectivity, proving the Lemma.
 \hspace{1 cm} $\blacksquare$

The proof of $\ast$), and thus of the Theorem, follows immediately
from Lemma 3 and Lemma 6.

\hspace{10cm}Q.E.D.

\section*{Reference}

\begin{enumerate}
\item J. Simons and D. Sullivan.  ``Structured Vector Bundles
Define Differential $K$-Theory''.  Quanta of Maths.  AMS and Clay
Mathematics Institute.  2010.  pp. 577-597.

\end{enumerate}

\end{document}